\input amstex
\input xy
\xyoption{all}
\documentstyle{amsppt}
\document
\magnification=1200
\NoBlackBoxes
\nologo
\hoffset1.5cm
\voffset2cm
\pageheight {17cm}

\def\mathcal{\Cal}
\def\mathbb{\bold}


\bigskip

\centerline{\bf On the derived category of $\overline{M}_{0,n}$}

\medskip

\centerline{Yu.~Manin, M.~Smirnov}
\medskip

\centerline{\it Max--Planck--Institut f\"ur Mathematik, Bonn, Germany}
\bigskip
{\it ABSTRACT.}  Using Keel's presentation and Orlov's theorem,
we give an inductive description of the derived category of
moduli spaces of $n$--pointed stable curves of genus zero
and some full exceptional collections in it. The detailed
calculations are given for $\overline{M}_{0,6}$.

\bigskip

\centerline{\bf  0. Introduction}

\medskip

{\bf 0.1. Summary.}  Let $f:\,X\to Y$ be a  monoidal transformation of smooth projective schemes
with smooth center $Z\subset Y$. Then $D^b(Coh\, X)$ in a definite sense can be described
in terms of   $D^b(Coh\, Y)$ and  $D^b(Coh\,Z)$: cf. [Or1]. 

\smallskip

We apply this strategy to the calculation of the derived categories of moduli spaces $\overline{M}_{0,n}$.
Since there are several methods of representing  $\overline{M}_{0,n}$ as the result of  a
sequence of blow--ups with controlled bases and centers (see e.g. [Ke] and [Ka]), one may later ask about 
the interaction of the results obtained etc.

\smallskip

In this note we focus on the representation discovered and used by Keel in [Ke]. Its
remarkable feature is that $\overline{M}_{0,n+1}$ is obtained from 
$\overline{M}_{0,n}\times \overline{M}_{0,4}$ by iterated blow--ups of 
manifolds isomorphic to  $\overline{M}_{0,p+1}\times \overline{M}_{0,q+1}$, $p+q=n$.
This produces a straightforward inductive description if one  
complements the tools used in the computation a version of  ``K\"unneth formula''
for the derived categories. The simplest useful form is furnished by
Proposition 2.1.18 in [B\"o] constructing exceptional collections on $V\times W$
from those on $V$ and $W$. A much more sophisticated version
is contained in [BoLaLu]; cf. also [Kel]. 

\smallskip

This approach might eventually
furnish a concise and functorial description of quivers
useful for computations of (and in) $D^b(Coh\,\overline{M}_{0,n}),$
but we could not achieve this goal so far.

\smallskip

We imagine two directions of future work on this subject.

\smallskip

(A) We might expect that replacing motives of $\overline{M}_{0,S}$ by their
derived (or enhanced derived) categories, we can get a better understanding of
their quantum cohomology: cf. [MaS].

\smallskip
(B) The moduli spaces  $\overline{M}_{0,S}$ are rigid as objects of commutative
geometry. Ample experience, including Gauss' $q$--numbers and quantum groups,
 shows that non--commutative (``quantum'')  deformations
of rigid objects are especially interesting. The beautiful recent paper [KeTe]
shows how to represent  $\overline{M}_{0,S}$ as the projective spectrum
of a Koszul graded ring of sections of its $log$--canonical bundle.
It would be interesting to study the deformations of the derived category
of  $\overline{M}_{0,S}$ in terms of the deformations of this coordinate ring,
as well as in terms of Bondal's quivers [Bo].
\smallskip

 See  [AuKaOr]  for a detailed discussion
of the case of del Pezzo surfaces, containing in particular $\overline{M}_{0,5}$.
\smallskip
Noncommutative deformations of general
$\overline{M}_{0,S}$, together with structural correspondences between them, that we envision, 
might become a sophisticated counterpart of the theory of 
 quantum groups.

\medskip

{\bf 0.1. Plan of the paper.} The first section summarizes the main theorem due to Orlov  ([Or1])
allowing one to reconstruct the derived category of a blow--up. We supplement it by some
computations in the derived category, stressing, in particular, the codimension
two blow--ups, that will appear in the applications to  $\overline{M}_{0,S}$.
The second section describes Keel's tower from [Ke]. In the third section the inductive
construction based on the results of Orlov and Keel is elaborated. It is complemented
by a new result (Proposition 3.5 and Corollary 3.6) allowing us in certain circumstances to get
full exceptional collections consisting  of locally free sheaves. Finally, sec. 4 
gives an application of this construction to  $\overline{M}_{0,6}$.

\medskip
{\it Acknowledgements.} We are very grateful to D.~Orlov for detailed and
stimulating answers to the questions of Yu.~M. Our thanks are due to Ch.~B\"ohning
for a discussion of ``Kunneth formula'' for derived categories.
\bigskip

\centerline{\bf 1. Orlov's theorem}

\medskip

{\bf 1.1. Exceptional collections.}  Let $K$ be a ground field;  $D$ a $K$--linear triangulated category;
$T$ its shift functor;
$E,F\in Ob\,D.$ Following [Bo], we  consider the graded complex of $K$--vector
spaces with trivial differential
$$
\roman{Hom}_D^{\bullet} (E,F):= \bigoplus_{k\in\bold{Z}} \roman{Hom}_D^{k} (E,F)[-k],\quad
\roman{Hom}_D^{k} (E,F):=\roman{Hom}_D (E,T^kF).
$$
Remind the following definitions.
\smallskip
(i) An object $E$ is called {\it exceptional one}, if
$\roman{Hom}_D^{\bullet} (E,E)$ consists of $K\cdot\roman{id}_E$ put in degree zero.
\smallskip

(ii) A family of exceptional objects $(E_{\alpha})$ indexed by a totally ordered
set of subscripts $\alpha$ is called
{\it an exceptional collection} (or an exceptional sequence), if
$$
\roman{Hom}_D^{\bullet} (E_{\alpha},E_{\beta})=0 \quad \roman{for}\quad  \alpha >\beta .
\eqno(1.1)
$$
It is called {\it full}, if it generates $D$.

\smallskip

(iii) An exceptional sequence  $(E_{\alpha})$ as above is called {\it strong},
if  for any $\alpha <\beta$ the complex  $\roman{Hom}_D^{\bullet} (E_{\alpha},E_{\beta})$
has at most one non--vanishing component which is then of degree zero.

\medskip

{\bf 1.2. Blow--ups.} Let $X$ be a smooth projective variety over $K$, $Y\subset X$ its
smooth closed subvariety, $j:\,Y\to X$ the respective  embedding. We consider
the diagram describing the monoidal transformation of $X$ with center $Y$:
$$
\xymatrix{\tilde{Y} \ar[r]^{i}     \ar[d]^{\pi}& 
 \tilde{X}\ar[d]^{q} \\
Y \ar[r]^{j}& X}
\eqno(1.2)
$$
Let $\Cal{N}$ be the  sheaf of sections of the normal bundle to $Y$ in $X$. 
The exceptional divisor $\widetilde{Y}$ in $\widetilde{X}$ is canonically isomorphic to 
the relative projective spectrum of the symmetric algebra 
$S_{\Cal{O}_Y}(\Cal{N}^t)$ of the dual sheaf of sections of $N$.
The rank $c\ge 2$ of $N$ equals to the codimension of $Y$;  fibers of $\pi$ are 
$\bold{P}^{c-1}$. This projective bundle carries the standard invertible relative sheaves 
$\Cal{O}_{\pi}(l)$.
\smallskip

From the general definitions, we obtain  the canonical exact sequences
$$
0\to \Cal{N}^t\to j^*(\Omega^1_X)\to \Omega^1_Y\to 0
\eqno(1.3)
$$
and
$$
0\to \pi^*(\Omega^1_Y)\to \Omega^1_{\widetilde{Y}}\to  \Omega^1_{\widetilde{Y}/Y}\to 0.
\eqno(1.4)
$$
\smallskip
The dual exact sequence to that in [Huy], p. 252, reads
$$
0\to  \Omega^1_{\widetilde{Y}/Y}(1)\to \pi^*(\Cal{N}^t) \to \Cal{O}_{\pi}(1)\to 0.
\eqno(1.5)
$$

\medskip

{\bf 1.3. From exceptional collections on $X$ and $Y$ to those on $\widetilde{Y}$
and $\widetilde{X}$.}  Generally, for a variety $Z$ we denote by $D(Z)$ the bounded derived category
$D^b(Coh\,Z)$.
\smallskip

Let $(E_{\alpha}), \alpha \in A,$ be an exceptional collection in $D(X)$ and $(F_{\beta}), \beta\in B,$ be an exceptional sequence
in $D(Y)$. Consider the following sequence of objects in $D(\widetilde{X})$,
that is indexed by the set of pairs $(\beta,l)$ for $\beta\in B, -c+1\le l\le -1$
and pairs $(\alpha ,l=0)$ for $\alpha\in A$:
$$
Ri_*(L\pi^*(F_{\beta})\otimes\Cal{O}_{\pi}(-c+1)),\dots , 
Ri_*(L\pi^*(F_{\beta})\otimes\Cal{O}_{\pi}(-1)),
Lq^*(E_{\alpha}).
\eqno(1.6)
$$
It is ordered as suggested in (1.6): in each group  with fixed $l$ according to the order
of $A$ or $B$, and groups are ordered by increasing $l$. 

\medskip

{\bf 1.4. Proposition.} {\it (i) The sequence (1.6) is exceptional.

\smallskip

(ii) If  $(E_{\alpha})$ is full in $D(X)$ and    $(F_{\beta})$ is full in $D(Y)$,
then (1.6) is full in $D(\widetilde{X})$.}

\smallskip

This is the Corollary 4.4 in [Or1] (notice that our notation for the blow--up
diagram (1.2) differs from Orlov's notation and is closer to the
Huybrechts notation in [Huy], Ch. 11). 

\smallskip

In the remaining part of this section we produce some formulas
for morphism spaces between remaining pairs of objects in (1.6).

\medskip

{\bf 1.5. Auxiliary results.} {\it (a) Adjunction/duality formula.} For a morphism $f:\,U\to V$ of smooth varieties, define
$\roman{dim}\,f:=\roman{dim}\,U - \roman{dim}\,V$ and put
$$
\omega_f:= \omega_U\otimes f^*(\omega_V^{-1}).
\eqno(1.7)
$$
Then for $F\in D(U), E\in D(V)$ we have functorial isomorphisms ([Huy], p. 87)
$$
\roman{Hom}_{D(V)}(Rf_*(F),E)\cong  \roman{Hom}_{D(U)}(F,Lf^*(E)\otimes \omega_f[\roman{dim}\,f])
\eqno(1.8)
$$
(tensor products by an invertible or more general locally free sheaf here and below need not be derived).

\smallskip

{\it (b) The sheaf $\omega_{\widetilde{Y}}$.}  We will prove that
$$
\omega_{\widetilde{Y}}\cong \pi^*\circ j^*(\omega_X)\otimes  \Cal{O}_{\pi}(-c).
\eqno(1.9)
$$
In fact, from (1.4) and (1.3) we get 
$$
\omega_{\widetilde{Y}}\cong \pi^*( j^*(\omega_X)\otimes \roman{det}\,\Cal{N})\otimes 
\omega_{\widetilde{Y}/Y}.
\eqno(1.10)
$$
From (1.5)  we obtain 
$$
\pi^*(\roman{det}\,\Cal{N})\cong  \omega_{\widetilde{Y}/Y}^{-1}\otimes \Cal{O}_{\pi}(-c).
\eqno(1.11)
$$
Substituting (1.11) to (1.10), we get (1.9).
\smallskip

{\it (c) The sheaf $\omega_i$. }  We wish to prove that
$$
\omega_i\cong\Cal{O}_{\pi}(-1),\quad  \omega_i[\roman{dim}\, i]\cong   \Cal{O}_{\pi}(-1)[-1].
\eqno(1.12)
$$
In fact, in the notation of diagram (1.2), we have
$$
\omega_{\widetilde{X}}\cong q^*(\omega_X)\otimes \Cal{O}_{\widetilde{X}}((c-1)\widetilde{Y}),\quad
i^*(\Cal{O}_{\widetilde{X}}(\widetilde{Y}))\cong \Cal{O}_{\pi}(-1),
\eqno(1.13)
$$
where we write $\widetilde{Y}$ in place of  the exceptional divisor $i(\widetilde{Y})$: cf. [Huy], p. 252.
\smallskip
From (1.7) and  (1.13) we find that
$$
\omega_i= \omega_{\widetilde{Y}}\otimes i^*(\omega_{\widetilde{X}}^{-1})\cong 
\omega_{\widetilde{Y}}\otimes i^*\circ q^*(\omega_X^{-1})\otimes \Cal{O}_{\pi}(c-1).
\eqno(1.14)
$$
Since $i^*\circ q^*= \pi^*\circ j^*$, substituting (1.9) into (1.14), we find (1.12).

\smallskip

{\it (d) The sheaf $\omega_{\pi}$. } From (1.11) we can calculate $\omega_{\pi}:=\omega_{\tilde{Y}}\otimes \pi^*(\omega_Y^{-1})\cong \omega_{\widetilde{Y}/Y}$:
$$
\omega_{\pi} \cong \roman{det}\,\pi^*(\Cal{N})^{-1}\otimes \Cal{O}_{\pi}(-c)
$$
and
$$
\omega_{\pi} [\roman{dim}\,\pi ]\cong \pi^*( \roman{det}\,\Cal{N})^{-1}\otimes \Cal{O}_{\pi}(-c) [c-1].
\eqno(1.15)
$$

\smallskip
{\it (e) The sheaf $\omega_j$.} We have
$$
\omega_j \cong \roman{det}\,\Cal{N}, \quad  \omega_j  [\roman{dim}\, j] \cong \roman{det}\,\Cal{N} [-c].
\eqno(1.16)
$$
This follows directly from (1.3).

\medskip

{\it (f) The sheaf $\omega_q$.}  Finally, from (1.7) and (1.13) we get
$$
\omega_q\cong  \Cal{O}_{\widetilde{X}}((c-1)\widetilde{Y})\cong  \omega_q[\roman{dim}\,q].
$$

\medskip

{\bf 1.6. Calculations of $\roman{Hom}$'s.} Below we will need some formulas for  the morphism
spaces between objects of the sequence (1.6) that do not follow directly from the exceptionality
of this sequence, i.~e.
$\roman{Hom}_{D(\widetilde{X})}^{\bullet}(G,H)$
where $G,H$ are objects of (1.6), indexed respectively by $(\beta_1 ,l_1)$ 
and $(\beta_2, l_2)$ with $l_1<l_2$, or else by $(\beta ,l), l<0,$ and $(\alpha ,0)$.
In fact, as Orlov has proved, for other pairs the respective $\roman{Hom}$'s come from
$Y, X$ or vanish.

\smallskip

Moreover, in our applications to Keel's tower, all centers of consecutive blow--ups have codimension
$c=2$. In this case the only relevant value of $l$ is $l=-1$ so that we need to consider
only the second option.
\smallskip
Here we take  two arbitrary objects $E\in D(X)$ and
$F\in D(Y)$. 

\smallskip

{\bf 1.7. Proposition.} {\it There are functorial in $E,F$ isomorphisms
$$
\roman{Hom}_{D(\widetilde{X})}(Ri_*(L\pi^*(F)\otimes \Cal{O}_{\pi}(l)),Lq^*(E))\cong 
\eqno(1.17)
$$
$$
\roman{Hom}_{D({Y})}(
F \otimes S_{\Cal{O}_Y}^{-l-1}(\Cal{N})) [1], Lj^*(E))\cong
\eqno(1.18)
$$
$$
\roman{Hom}_{D({X})}(Rj_*
(F \otimes S_{\Cal{O}_Y}^{-l-1}(\Cal{N})\otimes \roman{det}\,\Cal{N}) [1-c], E).
\eqno(1.19)
$$
}
\smallskip

{\bf Proof.} Applying (1.8) to $i:\,\widetilde{Y}\to \widetilde{X}$ in place of $f$, we see that (1.17) is iso\-morphic to
$$
\roman{Hom}_{D(\widetilde{Y})}(L\pi^*(F)\otimes \Cal{O}_{\pi}(l),Li^*\circ Lq^*(E)\otimes \omega_i[-1]).
\eqno(1.20)
$$
Replacing here $\omega_i$ by (1.12) and $Li^*\circ Lq^*$ by  $L\pi^*\circ Lj^*$, we rewrite (1.20) as
$$
\roman{Hom}_{D(\widetilde{Y})}(L\pi^*(F)\otimes \Cal{O}_{\pi}(l),L\pi^*\circ Lj^*(E)
\otimes \Cal{O}_{\pi}(-1)[-1]).
\eqno(1.21)
$$
Multiply both arguments of (1.21) by the same  invertible sheaf
$\pi^*(\roman{det}\,\Cal{N})^{-1}\otimes \Cal{O}_{\pi}(1-c)$ and then apply to them the same
shift  $[c]$, without changing $\roman{Hom}$. Using (1.15) to rewrite the new
second argument, we see that the result will be
$$
\roman{Hom}_{D(\widetilde{Y})}(L\pi^*(F)\otimes \pi^*(\roman{det}\,\Cal{N})^{-1}\otimes
\Cal{O}_{\pi}(l+1-c)[c],  L\pi^*\circ Lj^*(E)\otimes \omega_{\pi} [\roman{dim}\,\pi ])
\eqno(1.22)
$$
We can rewrite (1.22) using the adjunction formula (1.8) for $\pi$:
$$
\roman{Hom}_{D({Y})}(R\pi_*\left[L\pi^*(F \otimes (\roman{det}\,\Cal{N})^{-1})\otimes
\Cal{O}_{\pi}(l+1-c)\right] [c],  Lj^*(E))
\eqno(1.23)
$$
Now apply the projection formula for $\pi$ (cf. [Huy], (3.11) on p. 83)
to the first argument:
$$
R\pi_*\left[L\pi^*(F \otimes (\roman{det}\,\Cal{N})^{-1})\otimes
\Cal{O}_{\pi}(l+1-c)\right] [c]
\cong
$$
$$ 
F \otimes (\roman{det}\,\Cal{N})^{-1}\otimes^L R\pi_*( \Cal{O}_{\pi}(l+1-c)) [c]
$$
and  remark that the complex   $R\pi_*( \Cal{O}_{\pi}(l+1-c))$
for $1-c\le l \le -1$ is quasi--isomorphic to $R^{c-1}\pi_*( \Cal{O}_{\pi}(l+c-1))[1-c]$.  Hence (1.23) will become
$$
\roman{Hom}_{D({Y})}(
F \otimes (\roman{det}\,\Cal{N})^{-1}\otimes R^{c-1}\pi_*( \Cal{O}_{\pi}(l+1-c)) [1], Lj^*(E)).
\eqno(1.24)
$$
Finally, relative Serre's duality and (1.5) imply that in the considered range of $l$  we have
$$
R^{c-1}\pi_*( \Cal{O}_{\pi}(l+1-c))\cong S_{\Cal{O}_Y}^{-l-1}(\Cal{N})\otimes \roman{det}\,\Cal{N}.
$$

Thus, we proved (1.18).
\smallskip
Multiplying both arguments by $\roman{det}\,\Cal{N}$, shifting by $-c$ and applying (1.16),
we get (1.19).

 This finishes the proof.

\medskip

{\it Remarks.} (i) For $l=-1,\, c=2$, (1.18) becomes simply
$$
\roman{Hom}_{D({Y})}(F [1], Lj^*(E)).
\eqno(1.25)
$$
This is the only case that must be considered, when $c=2$. We will use this formula in sec. 3 and 4.

\medskip

(ii) We could have started a proof of Prop. 1.7 by applying formula (1.8) to $q$, rather than to $i$,
thus following the diagram (1.8) clockwise rather than counter--clockwise.

\medskip

{\bf 1.8. Reconstruction of the derived category of a blow--up.} For simpli\-city, we continue 
considering the case of codimension 2 blow--up. The main result of [Or1] in this case
establishes  the canonical semiorthogonal composition of $D(X)$, whose components are
admissible subcategories:
$(\widetilde{D}(Y)_{-1}, D(X)_0)$, where $D(X)_0:=Rq^*(D(X))$ and 
$\widetilde{D}(Y)_{-1}$ consists  of objects 
$Ri_*(L\pi^*(F)\otimes\Cal{O}_{\pi}(-1))$ where $F\in Ob\,D(Y)$.
\smallskip
Together with the formula (1.17)$\cong$(1.19) for $c=2, l=-1$, this semiorthogonal decomposition allows us
to reconstruct $D(\widetilde{X})$, however, not directly, but only
through the mediation of {\it DG--categories}, say, of the canonical enhancements
of all relevant derived categories (cf. [BoKa]) and [LuOr]).

\smallskip

The relevant construction on the level of DG--categories is described in [Ta], Definition 3.57,
under the name ``upper triangular DG--categories''.
Here we briefly describe its relation to the (two terms) semiorthogonal decompositions
of triangular categories, as was explained to us in [Or2].

\smallskip

Generally, let $\bold{A}$ be a linear pretriangulated DG--category, and $\Cal{A}:=H^0(\bold{A})$
its homotopy category with its natural triangulated structure. Assume that
$\Cal{A}$ admits a semiorthogonal decomposition $(\Cal{E},\Cal{F})$ as above, in the case of blow--up.
Lift $\Cal{F},\Cal{E}$ to the subcategories $\bold{F},\bold{E}$ of $\bold{A}$,
and consider a pair of objects $F,E$ of the respective categories.
The bifunctor given on objects by  $(F,E)\mapsto Hom_{\bold{A}}(F,E)$ with values in complexes of linear spaces
can be considered as DG--functor $\bold{\Phi}:\, \bold{E}\to \bold{Mod}$--$\bold{F}$,
``the upper right corner of the upper triangular DG--category'' in the sense of [Ta].

\smallskip

Conversely, given a triple $(\bold{F},\bold{E},\bold{\Phi})$ as above,
one can explicitly construct a new DG--category  $\bold{B}$. Its objects are triples
$(F,E,q)$ where $q\in \bold{\Phi}(E)(F)^1$. Moreover, morphisms constitute the complex
$$
Hom_{\bold{B}}((F,E,q), (F^{\prime},E^{\prime}, q^{\prime})) :=
Hom_{\bold{F}}(F,F^{\prime})\oplus Hom_{\bold{E}}(E,E^{\prime})\oplus 
\bold{\Phi}(E^{\prime},F)
$$
with pretty obvious composition rule.

\smallskip

Now, the principle result is that {\it $\bold{B}$ is quasiequivalent to $\bold{A}$.}

\smallskip

This fact, together with uniqueness of enhancements and explicit description
(1.25) of the bifunctor $\bold{\Phi}$ in the blow--up case, furnishes
the reconstruction of the triangulated category.

\bigskip
\centerline{\bf  2. Keel's tower}

\medskip

{ \bf 2.1. Notation: combinatorics of marks.} Let $S$ be a finite set,  $\roman{card}\,S=n\ge 3$.
We will call {\it an inductive structure on $S$} the choice of a three--element
subset $P\subset S$. We will sometimes denote $\Sigma:=S\setminus P$
so that $S=\Sigma\sqcup P$.
\smallskip

Recall that boundary strata of $\overline{M}_{0,S}$ are bijectively numbered
by the (isomorphism classes of)  stable $S$--marked trees. Such a tree
describes the dual combinatorial type of the curve parametrized by the generic
point of the respective stratum. The number of edges of such a tree 
equals the codimension of the stratum.
\smallskip
In particular,
boundary divisors, that is, trees with one edge, are determined by {\it
unordered 2--partitions $S=S_1\sqcup S_2$}, stable in the sense that
$|S_i|\ge 2$. Furthermore, codimension two strata are determined
by {\it 3--partitions} $S=S_1\sqcup S_2\sqcup S_3$ in which the middle term $S_2$
is uniquely defined, whereas $S_1$ and $S_2$ can be interchanged.
Stability condition here means that $|S_1|, |S_3|\ge 2, |S_2|\ge 1.$
\smallskip

Whenever an inductive structure $P$ is chosen on $S$, and $|S|\ge 4$, we may and will
{\it order each stable 2--partition} by the condition $|S_1\cap P|\le 1$, and for $|S|\ge 4$
we will {\it order  each stable  3--partition} by the condition $|(S_1\cup S_2)\cap P|\le 1.$

\medskip

{\bf 2.2. Keel's blow--down.} Now, assuming $n:=|S|\ge 4$ and given an inductive structure
on $S$, consider  a one--point set $\{\bullet\}$ disjoint from $S$
and the diagram of two  forgetful morphisms, forgetting  respectively
sections marked by $\Sigma$ and the one marked by $\bullet$:
$$
\xymatrix{ \overline{M}_{0,S\sqcup \{\bullet\}}\ar[r]^{f_{\Sigma}} \ar[d]^{f_{\{\bullet\}}} 
& \overline{M}_{0,P\sqcup \{\bullet\}}\\
 \overline{M}_{0,\Sigma\sqcup P}  }
$$

Notice that  $\overline{M}_{0,P\sqcup \{\bullet\}}$ is $\bold{P}^1$ endowed with three boundary
points. They are canonically marked by unordered partitions of $P\sqcup \{\bullet\}$
into two parts of cardinality 2. Such a partition, in turn, is determined by an  element $p$ of
$P$ (one part is $ \{\bullet ,p\}$) or a two--element subset of $P$ (the part not containing $\bullet$).

\smallskip
We summarize below some results of [Ke], showing, in particular,
that the map  (inverse to)
$$
(  f_{\{\bullet\}}, f_{\Sigma}):\,   \overline{M}_{0,S\sqcup\{\bullet\}}\to
\overline{M}_{0,S}\times
\overline{M}_{0,P\sqcup \{\bullet\}}
  \eqno(2.1)
$$
is a composition of blow--ups of smooth codimension two subvarieties
isomorphic to boundary divisors of $\overline{M}_{0,S}$.
These blow ups  naturally form a sequence of $n-3$ steps $b_k:\,B_{k,S}\to B_{k-1,S}$,
$k=2,\dots , n-2$. At each step, a union of pairwise disjoint connected smooth
submanifolds of codimension two is blown up. 

\smallskip

In order to bridge our notation with Keel's, the reader should have in mind the following
case:
$$
S:=\{1,\dots ,n\},\quad \bullet:= n+1, \quad  P:=\{1,2,3\},\quad  \Sigma := \{4,\dots ,n\} .
\eqno(2.2)
$$
\medskip
{\bf 2.3. Exceptional divisors of $( f_{\{\bullet\}}, f_{\Sigma})$.} In our notation, 
Lemma 1 of [Ke] establishes that exceptional divisors of the morphism
$( f_{\{\bullet\}}, f_{\Sigma})$ are exactly all boundary divisors of
$\overline{M}_{0,S\sqcup\{\bullet\}}$ corresponding to the stable 
2--partitions of  $S\sqcup \{\bullet\}$, satisfying the following condition:

\smallskip

{\it (*)  The part containing $\bullet$ contains no more than one element of $P$
and has cardinality $\ge 3$.}

\smallskip

As an independent check, the reader can convince oneself that
the number of such partitions coincides with the difference of
ranks of the Picard groups
$$
\roman{rk\,Pic\,} \overline{M}_{0,S\sqcup\{\bullet\}} -\roman{rk\,Pic\,}
(\overline{M}_{0,S}\times
\overline{M}_{0,P\sqcup \{\bullet\}})= 2^{n-1}-n-1.
$$

Let $\overline{\sigma}$ be a partition of $S\sqcup \{\bullet\}$ satisfying (*) above,
and let $\sigma$ be be the respective partition of $S$ obtained by deleting $\bullet$.
Obviously, we have
$$
f_{\bullet} (D_{\overline{\sigma}})=D_{\sigma}.
\eqno(2.3)
$$
We will call the cardinality of the second part of $\overline{\sigma}$
(and of $\sigma$) {\it the height} of $D_{\overline{\sigma}}$.
\medskip

{\bf 2.4. Keel's tower.} The main result of [Ke], sec. 1,  can now be stated in the following way.

\smallskip

The morphism (2.1) can be represented as a product of blowing--downs
$$
 \overline{M}_{0,S\sqcup\{\bullet\}}=:B_{n-2,S}\to B_{n-3,S}\to \dots \to B_{1,S}:=
\overline{M}_{0,S}\times\overline{M}_{0,P\sqcup \{\bullet\}}
\eqno(2.4)
$$
satisfying the following conditions: 

\smallskip

(i) Image of any exceptional divisor $D_{\overline{\sigma}}$ of height $h$ remains a divisor in  $B_{h,S}$,
but becomes a closed subscheme of codimension 2 in  $B_{h-1,S},\dots ,B_{1,S}$.
The product of the subsequent arrows, followed by the projection $B_{1,S}\to \overline{M}_{0,S}$
identifies this subscheme with $D_{\sigma}$.

\smallskip

(ii) Each morphism $B_{h+1,S}\to B_{h,S}$ is the blow up of the disjoint union
of those subschemes in $B_{h,S}$ that are images of exceptional divisors of height $h+1$.
Hence connected components of the center of the respective blow up
are isomorphic to $\overline{M}_{0,p+1}\times \overline{M}_{0,q+1}$, $p,q\le n-2$.

\bigskip

\centerline{\bf 3. Inductive construction of semiorthogonal decompositions}
\smallskip
\centerline{\bf and  complete 
exceptional collections.} 

\medskip

{\bf 3.1. The inductive step I: functoriality in $S$.}  In order to calculate $D^b(\overline{M}_{0,S\sqcup\{\bullet\}})$
assuming that the derived categories of the respective
modular spaces of smaller dimension are already known, we will apply
D.~Orlov's results summarized in sec.~1 to Keel's tower.
\smallskip
More precisely, Keel's tower depends  on the choice of an {\it inductive structure}
on $S$ in the sense of 2.1. In order to be able to induce a given inductive structure
on the subsets of $S$, we will adopt here the following convention, essentially returning us to the Keel's choice
(2.2). 

\smallskip

{\it $S$ is totally ordered, and for $|S|\ge 3$, $P\subset S$ consists of the first three elements of $S$.}
\smallskip
The inductive structure induced on  subsets of $S$ is defined then by the induced order.

\smallskip

With this conditions, one easily sees that a bijection of two sets of marks
compatible with their respective orders lifts to a unique isomorphism of Keel's towers.

\medskip

{\bf 3.2. The inductive step II: Keel's blow--ups.} Our ``inductive leap'' from
$S$ to $S\sqcup\{\bullet\}$ breaks down into the sequence of small inductive steps
corresponding to the consecutive floors of the Keel tower (2.4).
They will allow us to obtain an inductive description of a class of semiorthogonal decompositions
of $D^b(Coh\,\overline{M}_{0,S})$. Concretely, in view of Orlov's theorem, each floor  $p_k:\,B_{k+1,S}\to B_{k,S}$
provides a semiorthogonal decomposition of $D^b(Coh\,B_{k+1,S})$
of the form
$$
D^b(Coh\,B_{k+1,S}):\quad (\widetilde{D}({Y}_{k,S})_{-1}, D(B_{k,S})_0)
\eqno(3.1)
$$
To be more precise, we recapitulate
Orlov's results in this context. Consider the following specializations of
the blow--up diagram (1.2):
$$
\xymatrix{\hat{Y}_{k,S} \ar[r]^{i_k}     \ar[d]^{\pi_k}& 
 B_{k+1,S}\ar[d]^{q_k} \\
Y_{k,S} \ar[r]^{j_k}& B_{k,S}}
\eqno(3.2)
$$
in which $i_k:\,\widehat{Y}_{k,S}
\to B_{k+1,S}$ is the embedding of the exceptional divisor of $q_k$. Since $q_k$ blows up
 codimension two submanifolds,  $\pi_k:\,\widehat{Y}_{k,S}\to Y_{k,S}$ is a
fibration with fiber $\bold{P}^1$. Let $O(-1)$ be the respective relative sheaf on  $\widehat{Y}_{k,S}$.
Then $D^b(Coh\,\widehat{Y}_{k,S})$ has the semiorthogonal decomposition 
$$
D^b(Coh\,\widehat{Y}_{k,S}):\quad (\, \pi_k^*(D^b(Coh\, {Y}_{k,S}))\otimes O(-1),  q_k^*(D^b(Coh\, {B}_{k,S}))\,)
\eqno(3.3)
$$

Now, the subcategory $\widetilde{D}({Y}_{k,S})_{-1}$ in (3.1) is defined as
$$
\widetilde{D}({Y}_{k,S})_{-1}:= Ri_{k*}[ \pi_k^*(D^b(Coh\, {Y}_{k,S}))\otimes O(-1)]
\eqno(3.4)
$$
whereas $D(B_{k,S})_0$ in (3.1) is
$$
D(B_{k,S})_0:= Lq_k^*(D^b(Coh\,B_{k,S})).
\eqno(3.5)
$$
Moreover, since $q_k$ blows up a disjoint union of smooth submanifolds $Y_{\sigma}$,
$\sigma = (S_1, S_2)$, each of the subcategories
$\widetilde{D}({Y}_{k,S})_{-1}$ admits an {\it orthogonal} decomposition
$$
\widetilde{D}({Y}_{k,S})_{-1}:\quad ( \widetilde{D}({Y}_{\sigma})_{-1} \,|\, \roman{card}\,S_2=k+2)
\eqno(3.6)
$$

\smallskip

Finally, we have a canonical identification
$$
Y_{\sigma}= \overline{M}_{0,S_1\sqcup \{\bullet_{\sigma}\}}\times
 \overline{M}_{0,S_2\sqcup \{\bullet_{\sigma}\}}
 \eqno(3.7)
 $$
 where $\bullet_{\sigma}$ corresponds to the intersection point of two components.
 Therefore $D^b(Coh\,Y_{\sigma})$ is generated
 by an external product $\boxtimes$ of any two exceptional collections generating respectively
$\overline{M}_{0,S_1\sqcup \{\bullet_{\sigma}\}}$ and $\overline{M}_{0,S_2\sqcup \{\bullet_{\sigma}\}}$.
\smallskip

This ``K\"unneth formula'' for derived and DG--categories is known on various levels
of generality. Ch.~B\"ohning ([B\"o])  establishes it for collections of locally free sheaves;
the proof generalizes to complexes of locally free sheaves. In  [BoLaLu] 
the DG case  is treated.
\smallskip

The base $B_{1,S}$  has a similar decomposition  (cf. (2.4)).
Thus, Keel's tower  provides a tool to generate semi--orthogonal
decompositions and exceptional collections for $D^b(\overline{M}_{0,n+1})$ 
from similar objects for  $D^b(\overline{M}_{0,m})$, $m\le n$.

\medskip

{\bf 3.3. Exceptional collections for small $n$.}  On $\overline{M}_{0,4}\equiv \bold{P}^1$
the standard full strong exceptional collection is $(\Cal{O}(-1), \Cal{O})$.

\smallskip

If we represent $\overline{M}_{0,5}$ as a blow--up $p:\,\overline{M}_{0,5} \to \bold{P}^2$ at four points,
and denote by $l_i$, $i=1,\dots ,4$ the respective exceptional divisors,
then for any choice $(F_0,F_1,F_2)$ of a full strong exceptional collection   on $\bold{P}^2$,
e.~g. $F_i=\Cal{O}(i-2)$,
Orlov's theorem will provide a full strong exceptional collection on $\overline{M}_{0,5}$ of the form
$$
( \Cal{O}_{l_1}(-1),\dots , \Cal{O}_{l_4}(-1), p^*F_0,p^*F_1,p^*F_2)
\eqno(3.8)
$$
In [KarN] it was shown that a sequence of mutations turns such a sequence into
a full exceptional collection of locally free sheaves (this trick
actually works for all del Pezzo surfaces).  Taking its $\boxtimes$ with the standard
collection on $\bold{P}^1$ and then using Keel's blow--up, one gets an
analog of (3.8) for  $\overline{M}_{0,6}$.

\smallskip

In the remaining part of this section we will show that under certain assumptions, one can find a sequence of mutations, turning the collection thus obtained into a collection
of  locally free sheaves. In the next section we will demonstrate that this procedure
can be applied  to $\overline{M}_{0,6}$.

\medskip

{\bf 3.4. Preparation.} Consider a blow-up diagram (1.2)
$$
\xymatrix{\tilde{Y}\ar[d]^{\pi} \ar[r]^{i}& 
 \tilde{X}\ar[d]^{q} \\
Y \ar[r]^{j}& X}
$$
For a general blow up $R^0q_*\Cal{O}_{\widetilde{X}}\simeq \Cal{O}_X$ and $R^iq_*\Cal{O}_{\widetilde{X}}=0$ (see [Or1], proof of Lemma 4.1). Hence the projection formula and the Leray spectral sequence give $H^i(\widetilde{X}, q^*\Cal{E})=H^i(X, \Cal{E})$ for any locally free sheaf $\Cal{E}$.
\smallskip

 From now on we assume that $j(Y)$ is of codimension two in $X$. If $F_1, \dots, F_r$
 (resp. $E_1, \dots, E_n$) is a full exceptional collection of locally free sheaves on $Y$ (resp. on $X$),
 then Orlov's collection (1.6)
$$
Ri_*(\pi^*F_1\otimes \Cal{O}_{\pi}(-1)), \dots, Ri_*(\pi^*F_r \otimes \Cal{O}_{\pi}(-1)), q^*E_1, \dots, q^*E_n
\eqno(3.9)
$$
is a full exceptional collection on $\widetilde{X}$.

\medskip

{\bf 3.5. Proposition.}   {\it Let $E_1, \dots E_n$ be a full exceptional collection of 
locally free sheaves on $X$. Assume moreover that $j^*E_1, \dots , j^*E_r$ form a full exceptional collection on $Y$. If $Y$ is of codimension $2$ in $X$, then the collection 
$$
q^*E_1,\, q^*E_1(\widetilde{Y}),\dots , q^*E_r \,, q^*E_r(\widetilde{Y}), q^*E_{r+1}, \dots, q^*E_n
\eqno(3.10)
$$
is a full exceptional collection on $\widetilde{X}$.}

\medskip

Here and below we use notation of the type $q^*L(\widetilde{Y})$ as
a shorthand for $q^*L\otimes \Cal{O}_{\widetilde{X}}(i(\widetilde{Y}))$.

\medskip

{\bf Proof.}  The strategy of our proof is simple. We start with the exceptional collection (3.9),
with $F_a:=q^*E_a$ for $1\le a\le r$, and show that it can be transformed into (3.10)
by an explicit sequence of mutations (cf. [Bo], [Kuz]).

\smallskip

More precisely, for $1\le a\le r$, put $A_a:=  Ri_*(\pi^*j^*E_a \otimes \Cal{O}_{\pi}(-1))$, $B_a:= q^*E_a.$

\smallskip

We will first check that
$$
\roman{Hom}^{\bullet}(A_b, B_a)=0\quad\roman{for}\quad b>a.
\eqno(3.11)
$$
so that the right mutation of such an exceptional pair simply reduces to the permutation  $(A_b,B_a)\mapsto (B_a, A_b).$
This  shows, that we may consecutively move $A_r$ in (3.9) to the right, until it reaches the position 
directly to the left of  $B_r$;  then  move $A_{r-1}$  to the right, until it reaches the position 
directly to the left of $B_{r-1}$; and so on. The  result will be the exceptional collection
$(A_1,B_1; A_2,B_2; \dots ; A_r,B_r; B_{r+1},\dots ,B_n)$.
 
 \smallskip
Second, we will check that 
$$
R_{B_a}(A_a)\cong B_a(\widetilde{Y})
\eqno(3.12)
$$
Thus additional $r$ right mutations will transform the latter collection into
$$
(B_1, B_1(\widetilde {Y})),  \dots ,  B_r, B_r(\widetilde {Y}), B_{r+1},\dots ,B_n)
$$ 
that is, to (3.10).

\medskip

{\it Proof of (3.11).}   Consider the isomorphism (1.17) $\cong$ (1.25) written
for $F:=j^*E_b$, $E:=E_a[i]$ where $i$ is an arbitrary shift and $b>a$. Its left hand side
will represent one of the components of $\roman{Hom}^{\bullet}(A_b, B_a)$. Hence it suffices to prove
that the right hand side vanishes. But it is simply $\roman{Hom}_{D(Y)}(j^*(E_b)[1], j^*E_a[i])$,
$1\le a\le b \le r$, and all these groups vanish, because we assumed that $j^*E_1,\dots , j^*E_r$
is an exceptional collection.

\medskip

{\it {Proof of (3.12).}} Consider now the case $a=b$. First of all, recall that 
$R_{B_a}(A_a)$ is defined as the cone $C(\alpha_a)$ of the canonical morphism in $D(\widetilde{X})$
$$
\alpha_a:\ A_a\to \roman{Hom}^{\bullet}(A_a,B_a)^t\otimes B_a
\eqno(3.13)
$$
where $t$ means linear dual in the category of  of graded linear spaces.
In order to calculate $\roman{Hom}^{\bullet}(A_a,B_a)$, we note that $A_a$ fits into exact sequence
$0\to  q^*E_a\to q^*E_a(\widetilde{Y})\to A_a\to 0$ (cf. (1.13)), and thus $A_a$
is quasi--isomorphic to its projective resolution
$$
\Cal{A}_a:\ 0\to q^*E_a\to q^*E_a(\widetilde{Y})\to 0
\eqno(3.14)
$$
(with the first term in degree $-1$). 

\smallskip

Since $A_a$ is exceptional, $\roman{Hom}^{\bullet}(\Cal{A}_a,B_a)$ is spanned
by the canonical isomorphism $\Cal{A}_a\to B_a[1]$ which is identity in degree $-1$,
and the morphism (3.13) can be represented by the morphism of complexes
$$
\xymatrix{
0\ar[d]\ar[r]&q^*E_a \ar[r]^{g}\ar[d]^{id}&q^*E_a(\widetilde{Y})\ar[d] \ar[r]&0\ar[d]\\
0\ar[r]&q^*E_a\ar[r]&0\ar[r]\ar[r]&0}
$$
The cone of it is the complex
$$
\xymatrix{
0\ar[r]&q^*E_a \ar[r]^<<<<<{(id, -g)}&q^*E_a \oplus q^*E_a(\widetilde{Y})\ar[r]&0},
$$
where $q^*E_a \oplus q^*E_a(\widetilde{Y})$ is in degree  $-1.$
There exists a short exact sequence of sheaves
$$
\xymatrix{
0\ar[r]&q^*E_a \ar[r]^<<<<<{(id, -g)}&q^*E_a \oplus q^*E_a(\widetilde{Y})\ar[r]^>>>>>{\psi}& q^*E_a(\widetilde{Y})\ar[r]&0, }
$$
where $\psi(v, w)=g(v)+w$. Therefore, the cone is quasi--isomorphic to the complex with one non--zero term $q^*E(\widetilde{Y})$ placed in degree $-1$, i.e. $q^*E_a(\widetilde{Y})[1]$. Hence $R_{B_a}(A_a)[1]\cong  q^*E_a(\widetilde{Y})[1]$ and finally
$R_{B_a}(A_a) \cong  B_a(\widetilde{Y}).$

\smallskip

This completes the proof of Proposition 3.5.

\bigskip

{\bf 3.6. Corollary.} {\it Let $E_1, \dots, E_n$ be a full exceptional collection of locally free sheaves on $X$. Assume moreover that
for some $s\le r$,  $j^*E_s, \dots , j^*E_r$ form a full exceptional collection on $Y$. If $Y$ is of codimension $2$ in $X$, then the collection
$$
q^*E_1(\widetilde{Y}), \dots, q^*E_{s-1}(\widetilde{Y}), q^*E_s,\, q^*E_s(\widetilde{Y}), \dots , q^*E_r \,, q^*E_r(\widetilde{Y}), q^*E_{r+1}, \dots, q^*E_n
$$
is a full exceptional collection on $\widetilde{X}$.}

\medskip

{\bf Proof.} First, let us recall a general fact. Let $\mathcal{D}$ be a triangulated category with a Serre functor $S\colon \mathcal{D} \to \mathcal{D} $ and $\mathcal{A}$ an admissible subcategory (cf. [BoKa1], [Kuz]). In this situation we have two semi-orthogonal decompositions $\langle \mathcal{A}, ^\perp\mathcal{A} \rangle$ and $\langle \mathcal{A}^\perp, \mathcal{A} \rangle$. By [BoKa1], Proposition 3.6, we obtain that $$
R_{^\perp\mathcal{A}}(\mathcal{A})=S^{-1}(\mathcal{A}),
\eqno(3.13)
$$
$$
L_{\mathcal{A}^\perp}(\mathcal{A})=S(\mathcal{A}).
\eqno(3.14)
$$
Our proof of corollary 3.6 will consist of applications of formulas (3.13),(3.14) and proposition 3.5.

\medskip
Let $\mathcal{D}=D(X)$, $\mathcal{A}=\langle E_1, \dots, E_{s-1} \rangle$, $^\perp\mathcal{A}=\langle E_s, \dots, E_{n} \rangle$ and the Serre functor is $S_X=\, \cdot \, \otimes \omega_X [\text{dim} X]$. Applying formula (3.13) we get a full exceptional collection
$$
E_s, \dots , E_n, S^{-1}_X(E_1), \dots , S^{-1}_X(E_{s-1}).
$$
To this collection we can apply proposition 3.5 and we get a full exceptional collection on $\widetilde{X}$
$$
q^*E_s, q^*E_s(\widetilde{Y}), \dots , q^*E_r, q^*E_r(\widetilde{Y}), q^*E_{r+1}, \dots , q^*E_n, Lq^* \circ S^{-1}_X(E_1), \dots , Lq^* \circ S^{-1}_X(E_{s-1}).
$$
Let now $\mathcal{D}=D(\widetilde{X})$, $\mathcal{A}=\langle Lq^*S^{-1}_X(E_1), \dots , Lq^*S^{-1}_X(E_{s-1}) \rangle$ and the Serre functor is $S_{\widetilde{X}}=\, \cdot \, \otimes \omega_{\widetilde{X}} [\text{dim} \widetilde{X}]$. Applying formula (3.14)  and using that
$$
S_{\widetilde{X}} \circ  Lq^* \circ S^{-1}_X \simeq \, Lq^* \otimes \mathcal{O}_{\widetilde{X}}(\widetilde{Y})
$$
we get the desired statement.

\newpage

\centerline{\bf 4. Example: moduli space $\overline{M}_{0,6}$}

\medskip

{\bf 4.1. Preparation: moduli space $\overline{M}_{0,5}$.} Let $S=\{1,2,3,4, 5\}$, $S'=\{1,2,3,4\}$ and $P=\{1,2,3\}$. Consider the Keel's tower
$$
\xymatrix{
B_{2,S'} = \overline{M}_{0,S} \ar[d]^{q_{1, S'}}\\
B_{1,S'} = \overline{M}_{0,S'}\times \overline{M}_{0,P\sqcup \{5\}}
}
$$
The map $q_{1, S'}$ contracts 3 boundary divisors $D_{\sigma_i}$ corresponding to unordered partitions
$$
\sigma_1=(5, 1, 4| 2 ,3) ,\quad \sigma_2=(5, 2, 4| 1 ,3) ,\quad \sigma_3=(5, 3, 4| 2 ,1).
$$

Let us identify   $\overline{M}_{0,S'}$ with $\overline{M}_{0,P\sqcup \{5\}}$ using the bijection of
labels identical on $\{1,2,3\}$ and mapping 4 to 5. Then images of $D_{\sigma_i}$ become three points
on the diagonal, corresponding to the partitions with one part $\{i,4\}$, resp. $\{i,5\}$.
 Let us imagine  $\overline{M}_{0,P\sqcup \{5\}}$  as the horizontal axis $\bold{P}^1$,
and  $\overline{M}_{0,S'}$ as the vertical one. Let $H_i$, resp. $V_i$, be the horizontal, resp. vertical fiber,
passing through the image of  $D_{\sigma_i}$.

\bigskip

Now denote by $\widetilde{H}_i$ and $\widetilde{V}_i$ the proper transforms of  $H_i$, resp. $V_i$,
in $ \overline{M}_{0,S}$, and let $\widetilde{Z}$ be the proper transform of the diagonal $Z$. Divisors $D_{\sigma_i}$, $\widetilde{H}_i$, $\widetilde{V}_i$ and $\widetilde{Z}$ are all isomorphic to $\mathbb{P}^1$ and have self--intersection $(-1)$. There are
precisely10 such curves on $\overline{M}_{0,S}$.

\medskip
Let $F_0, F_1$ and $G_0, G_1$ be full exceptional collections of locally free sheaves on $\overline{M}_{0,S'}$ and $\overline{M}_{0,P\sqcup \{5\}}$ respectively. It is well known that $F_i$, $G_i$ invertible sheaves. We know that
$$
F_0\boxtimes G_0, F_1\boxtimes G_0, F_0\boxtimes G_1, F_1\boxtimes G_1
$$
is a full exceptional collection on $\overline{M}_{0,S'}\times \overline{M}_{0,P\sqcup \{5\}}$. Denote it as
$$
L_0, L_1, L_2, L_3.
\eqno(4.1)
$$

Consider the decomposition $q_{1, S'}=f_1\circ f_2 \circ f_3$ where $f_i$ contracts only $D_{\sigma_i}$. Each $f_i$ is a blow-up of a surface at a point. Applying corollary 3.6 to the blow-up $f_3$ and using restriction of $L_0$ as a full exceptional collection in $D(pt)$ we obtain a full exceptional collection on the resulting surface
$$
f_3^*L_0, f_3^*L_0(D_{\sigma_3}), f_3^*L_1, f_3^*L_2, f_3^*L_3.
$$
Continuing in the same way and always using restriction of the first element as a full exceptional collection in $D(pt)$
we obtain a full exceptional collection on $\overline{M}_{0,S}$
$$
q_{1, S'}^*L_0, q_{1, S'}^*L_0(D_{\sigma_1}), q_{1, S'}^*L_0(D_{\sigma_2}), q_{1, S'}^*L_0(D_{\sigma_3}), q_{1, S'}^*L_1, q_{1, S'}^*L_2, q_{1, S'}^*L_3.
\eqno(4.2)
$$
Of course, this is just an example of a full exceptional collection on $\overline{M}_{0,S}$. If we used restrictions of other elements of (4.1) we would have obtained a different answer.

\medskip
{\bf 4.2. Preparation: moduli space $\overline{M}_{0,6}$.} Here we will introduce convenient notations which will be used later to write down a full exceptional collection on $\overline{M}_{0,6}$.

Let $S=\{1,2,3,4,5\}$, $\{\bullet\}=\{6\}$ and $P=\{1,2,3\}$. Consider the Keel's tower
$$
\xymatrix{
B_{3,S} = \overline{M}_{0,S\sqcup \{6\}} \ar[d]^{q_{2,S}}\\
B_{2,S}\ar[d]^{q_{1,S}} \\
B_{1,S} = \overline{M}_{0,S}\times \overline{M}_{0,P\sqcup \{6\}}
}
\eqno(4.3)
$$
The map $q_{1,S} \circ q_{2,S}$ contracts 10 boundary divisors $E_i$, $1 \leq i \leq 10$. At the height 3 level it contracts 7 boundary divisors corresponding to the following unordered partitions
$$
E_4 \leftrightarrow (6, 1, 4 | 5, 2 ,3), \quad E_6 \leftrightarrow (6, 2, 4 | 5, 1 ,3), \quad  E_8 \leftrightarrow (6, 3, 4 | 5, 2 ,1),
$$
$$
E_5 \leftrightarrow (6, 1, 5 | 4, 2 ,3), \quad  E_7 \leftrightarrow (6, 2, 5 | 4, 1 ,3), \quad  E_9 \leftrightarrow (6, 3, 5 | 4, 2 ,1),
$$
and $E_{10} \leftrightarrow (6, 4, 5 | 1, 2 ,3)$. At the height 2 level it contracts images under $q_{2,S}$ of 3 boundary divisors corresponding to the following unordered partitions
$$
E_1 \leftrightarrow (6, 1, 4, 5| 2 ,3), \quad E_2 \leftrightarrow (6, 2, 4, 5| 1 ,3), \quad E_3 \leftrightarrow (6, 3, 4, 5| 2 ,1).
$$

The divisors  $E_1, E_2, E_3$ are pairwise disjoint and $E_4,\dots E_{10}$ are pairwise disjoint as well. We list below all nonempty intersections
$$
E_1 \cdot E_4 =P_1,  \quad E_1 \cdot E_5=Q_1,
$$
$$ 
E_2 \cdot E_6 =P_2, \quad E_2 \cdot E_7=Q_2,
$$
$$
E_3 \cdot E_8 =P_3, \quad E_3 \cdot E_9=Q_3,
$$
$$
E_1 \cdot E_{10}=R_1; \quad E_2 \cdot E_{10}=R_2; \quad E_3 \cdot E_{10}=R_3,
$$
where $P_i$, $Q_i$, $R_i$ are isomorphic to $\mathbb{P}^1$ and all intersections are transversal.

\medskip
Let $L_0, \dots, L_6$ be an exceptional collection of invertible sheaves on $\overline{M}_{0,S}$ and $G_0, G_1$ an exceptional collection of invertible sheaves on $\overline{M}_{0,P\sqcup \{6\}}$. From them we construct a collection on $\overline{M}_{0,S}\times \overline{M}_{0,P\sqcup \{6\}}$ 
$$
L_0\boxtimes G_0, \dots L_6\boxtimes G_0, L_0\boxtimes G_1, \dots L_6\boxtimes G_1.
\eqno(4.4)
$$
We will need some assumptions on the collection $L_0, \dots , L_6$. We use notations for $\overline{M}_{0,S'}$ introduced earlier in section 4.1.

\medskip
{\it Assumption 1.} $L_1, L_2$ restricted to $D_{\sigma_1}$ form a full exceptional collection on $D_{\sigma_1}$; $L_2, L_3$ restricted to $D_{\sigma_2}$ form a full exceptional collection $D_{\sigma_2}$; $L_3, L_4$ restricted to $D_{\sigma_3}$ form a full exceptional collection $D_{\sigma_3}$.

\medskip

{\it Assumption 2.} $L_0, L_1$ restrict to a full exceptional collection on $\widetilde{H}_1$ and $\widetilde{V}_1$. The same holds for $L_1, L_2$ on $\widetilde{H}_2, \widetilde{V}_2$ and  $L_2, L_3$ on $\widetilde{H}_3, \widetilde{V}_3$ and $L_5, L_6$ on $\widetilde{Z}$.

\medskip
These assumptions are satisfied, for example, for
$$
\mathcal{O}, \mathcal{O}(D_{\sigma_1}), \mathcal{O}(D_{\sigma_2}), \mathcal{O}(D_{\sigma_3}), q_{1,S'}^*\mathcal{O}(0,1), q_{1,S'}^*\mathcal{O}(1,0), q_{1,S'}^*\mathcal{O}(1,1),
$$
where we used notations of section 4.1 and identification of $\overline{M}_{0,S'}\times \overline{M}_{0,P\sqcup \{5\}}$ with $\mathbb{P}^1\times \mathbb{P}^1$.

\medskip
{\bf 4.3. Collection.} Similar to section 4.1, in view of Keel's tower (4.3), consecutively applying corollary 3.6 to full exceptional collection (4.4) one can obtain a full exceptional collection on $\overline{M}_{0,6}$. We start with listing its elements,
and afterwards give some indications about how it was obtained.

\medskip
Let $q=q_{1,S} \circ q_{2,S}$ and $E_{\geq i}=\sum_{k=i}^{k=10}E_k.$
The exceptional collection is
$$
q^*(L_0\boxtimes G_0)(E_{\geq 1}),
$$
$$
q^*(L_1\boxtimes G_0)(E_{\geq 2}), q^*(L_1\boxtimes G_0)(E_{\geq 1}),
$$
$$
q^*(L_2\boxtimes G_0)(E_{\geq 2}), q^*(L_2\boxtimes G_0)(E_1+E_{\geq 3}), q^*(L_2\boxtimes G_0)(E_{\geq 1}),
$$
$$
q^*(L_3\boxtimes G_0)(E_{\geq 3}), q^*(L_3\boxtimes G_0)(E_2+E_{\geq 4}), q^*(L_3\boxtimes G_0)(E_{\geq 2}),
$$
$$
q^*(L_4\boxtimes G_0)(E_{\geq 4}), q^*(L_4\boxtimes G_0)(E_{\geq 3}),
$$
$$
q^*(L_5\boxtimes G_0)(E_{\geq 4}),
$$
$$
q^*(L_6\boxtimes G_0)(E_{\geq 4}), 
$$
\smallskip
$$
q^*(L_0 \boxtimes  G_1)(E_{\geq 5}), q^*(L_0 \boxtimes  G_1)(E_4+E_{\geq 6}), q^*(L_0 \boxtimes  G_1)(E_{\geq 4}), 
$$
$$
q^*(L_1 \boxtimes  G_1)(E_{\geq 6}), q^*(L_1 \boxtimes  G_1)(E_{\geq 5}), q^*(L_1 \boxtimes  G_1)(E_4+E_{\geq 7}), q^*(L_1 \boxtimes  G_1)(E_4+E_6+E_{\geq 8}), 
$$
\hfill{$q^*(L_1 \boxtimes  G_1)(E_4+E_{\geq 6})$,}
$$
q^*(L_2 \boxtimes  G_1)(E_{\geq 8}), q^*(L_2 \boxtimes  G_1)(E_{\geq 7}), q^*(L_2 \boxtimes  G_1)(E_6+E_{\geq 9}), q^*(L_2 \boxtimes  G_1)(E_6+E_8+E_{10}), 
$$
\hfill{$q^*(L_2 \boxtimes  G_1)(E_6+E_{\geq 8})$},
$$
q^*(L_3 \boxtimes  G_1)(E_{10}), q^*(L_3 \boxtimes  G_1)(E_{\geq 9}), q^*(L_3 \boxtimes  G_1)(E_8+E_{10}),
$$
$$
q^*(L_4 \boxtimes  G_1)(E_{10}),
$$
$$
q^*(L_5 \boxtimes  G_1),q^*(L_5 \boxtimes  G_1)(E_{10}), 
$$
$$
q^*(L_6 \boxtimes  G_1),q^*(L_6 \boxtimes  G_1)(E_{10}).
$$
Below we describe the algorithm that was used to obtain this collection.

\medskip
{\bf 4.3.1 Algorithm, step I.} Due to Assumption 1,  the restriction of the pair
$$
L_1 \boxtimes G_0, L_2 \boxtimes G_0
$$
to $q(E_1)$ is a full exceptional collection. The same holds for pairs $L_2 \boxtimes G_0, L_3 \boxtimes G_0$ on $q(E_2)$ and $L_3 \boxtimes G_0, L_4 \boxtimes G_0$ on $q(E_3)$.
\smallskip
Represent $q_{1,S}=f_1 \circ f_2 \circ f_3$ as a composition of three blow--ups (cf. ~section 4.1) in such a way that the (preimage of) $q(E_i)$ is blown up at the $i$--th step for $1 \leq i \leq 3$. At the first step we use the pair $L_1 \boxtimes G_0, L_2 \boxtimes G_0$ to apply Corollary 3.6. At the second step we use the pair
$$
f_1^*(L_2 \boxtimes G_0)(f_2\circ f_3\circ q_{2,S}(E_1)), f_1^*(L_3 \boxtimes G_0),
$$
which restricts to a full exceptional collection on $f_2\circ f_3\circ q_{2,S}(E_2)$ because $q(E_1)$ and $q(E_2)$ are disjoint and $L_2 \boxtimes G_0, L_3 \boxtimes G_0$ restricts to an exceptional collection on $q(E_2)$ as was pointed out above.
\smallskip

One proceeds similarly at the third step and obtains the following exceptional collection on $B_{2,S}$
$$
q_{1,S}^*(L_0\boxtimes G_0)(E'_1+E'_2+E'_3),
$$
$$
q_{1,S}^*(L_1\boxtimes G_0)(E'_2+E'_3), q_{1,S}^*(L_1\boxtimes G_0)(E'_1+E'_2+E'_3), 
$$
$$
q_{1,S}^*(L_2\boxtimes G_0)(E'_2+E'_3), q_{1,S}^*(L_2\boxtimes G_0)(E'_1+E'_3), q_{1,S}^*(L_2\boxtimes G_0)(E'_1+E'_2+E'_3), 
$$
$$
q_{1,S}^*(L_3\boxtimes G_0)(E'_3), q_{1,S}^*(L_3\boxtimes G_0)(E'_2), q_{1,S}^*(L_3\boxtimes G_0)(E'_2+E'_3), 
$$
$$
q_{1,S}^*(L_4\boxtimes G_0), q_{1,S}^*(L_4\boxtimes G_0)(E'_3), 
$$
$$
q_{1,S}^*(L_5\boxtimes G_0), q_{1,S}^*(L_6\boxtimes G_0), 
$$
$$
q_{1,S}^*(L_0 \boxtimes  G_1), q_{1,S}^*(L_1 \boxtimes  G_1), \dots,
$$
where $E_i'=q_{2,S}(E_i)$.

\medskip
{\bf 4.3.2 Algorithm, step II.} Due to Assumption 2, the restriction of the pair
$$
L_0 \boxtimes G_1, L_1 \boxtimes G_1
$$
to $q_{2,S}(E_4)$ and $q_{2,S}(E_5)$  gives full exceptional collections on them. The same holds for the pair $L_1 \boxtimes G_1, L_2 \boxtimes G_1$ on $q_{2,S}(E_6)$ and $q_{2,S}(E_7)$; for $L_2 \boxtimes G_1, L_3 \boxtimes G_1$ on $q_{2,S}(E_8)$ and $q_{2,S}(E_9)$; for $L_5 \boxtimes G_1, L_6 \boxtimes G_1$ on $q_{2,S}(E_{10})$.

\smallskip

Represent $q_{2,S}=g_4\circ \dots \circ g_{10}$ as the composition of  blow--downs $g_i$ of $E_i$ for  $4\leq i \leq 10$. To apply Corollary 3.6 to a single blow--up one needs to choose an exceptional pair. We always choose pairs related to those described in the beginning of this section(in fact, they are pull--backs of those followed by a twist with $\mathcal{O}(D)$, where $D$ is a divisor disjoint from the exceptional divisor considered at this step).

\bigskip

\centerline{\bf References}

\medskip

[AuKaOr] D.~Auroux, L.~Katzarkov, D.~Orlov. {\it Mirror
symmetry for del Pezzo surfaces: vanishing cycles and coherent sheaves.}
Invent. Math. 166(2006), no.  3, 537--582.
arXiv:math/0506166

\smallskip

[Bo] A.~Bondal. {\it Representations of associative algebras and coherent sheaves.}
Math. USSR Izv. 34 (1990), 23--42.

\smallskip

[BoKa1] A.~Bondal, M.~Kapranov. {\it  Representable functors, Serre functors, and reconstructions.} (Russian) Izv. Akad. Nauk SSSR Ser. Mat. 53 (1989), no. 6, 1183--1205; translation in Math. USSR Izv. 35 (1990), no. 3, 519--541

\smallskip

[BoKa2] A.~Bondal, M.~Kapranov. {\it Framed triangulated categories.} Math USSR Sbornik, 70, Nr 1 (1991), 93--107.

\smallskip

[BoLaLu]  A.~Bondal, M.~Larsen, V.~Lunts. {\it Grothendieck ring of pretriangulated categories.}
Int.~Math.~Res.~Notes, 29 (2004), 1461--1495.
\smallskip
[B\"o] Ch.~B\"ohning. {\it Derived categories of coherent sheaves on rational homogeneous manifolds.}
Documenta Mathematica, 11 (2006), 261--331.

\smallskip
[Huy] D.~Huybrechts. {\it Fourier--Mukai transforms in algebraic geometry.} Oxford UP, 2006.

\smallskip

[Ka] M.~Kapranov.  {\it Veronese curves and Grothendieck--Knudsen moduli space $\overline M_{0,n}$.}
  J. Algebraic Geom.  2  (1993),  No. 2, 239--262.

\smallskip

[KarN] B.~V.~Karpov, D.~Yu.~Nogin. {\it Three--block exceptional collections over
del Pezzo surfaces.} Izv. Math. 62 (1998), no. 3, 429Ð-463.  arXiv:alg-geom/9703027

\smallskip

[Ke] S.~Keel. {\it Intersection theory of moduli space of stable $N$--pointed  curves of genus zero.}
Trans. AMS, Vol. 330, No. 2 (1992), 545--574.

\smallskip

[KeTe] S.~Keel, J.~Tevelev. {\it Equations for $\overline{M}_{0,n}.$} Int. Journ. Math.,
Vol. 20, No. 9 (2009), 1159--1184.

\smallskip

[Kel] B.~Keller.  {\it  On differential graded categories.}
International Congress of Mathematicians. Vol. II, 151--190, Eur. Math. Soc., ZŸrich, 2006. arXiv:0601185

\smallskip

[Kuz] A.~Kuznetsov. {\it Derived categories of cubic fourfolds.} In:
Cohomological and geometric approaches to rationality problems, 
Progr. Math., 282, BirkhŠuser Boston, Inc., Boston, MA, 2010, 219--243, arXiv:0808.3351

\smallskip

[LuOr] V.~Lunts, D.~Orlov. {\it Uniqueness of enhancement for triangulated categories.}
Journ. of the AMS, vol. 23, No 3 (2010), 853--908.

\smallskip

[MaS]  Yu.~Manin, M.~Smirnov.  {\it Towards motivic quantum cohomology of $\overline{M}_{0,S}$.}
arXiv:1107.4915

\smallskip

[Or1]  D.~Orlov. {\it Projective bundles, monoidal transformations, and derived categories of coherent sheaves.}
Russian Acad. Sci. Izv. Math., vol. 41 (1993), No. 1,  133--141.

\smallskip

[Or2] D.~Orlov. {\it Letter to Yu. Manin of 19/09/2011.}

\smallskip

[Ta] G.~Tabuada. {\it Th\'eorie homotopique des DG--cat\'egories.}  arXiv:0710.4303

\enddocument